\documentclass[12pt]{article}

\usepackage{amssymb}
\usepackage{amsmath}

\catcode`@=11
\long\def\@makefntext#1{\noindent #1}
\def\@evenfoot{}\def\@oddfoot{}
\def\@evenhead{\hbox to\textwidth{\footnotesize\rm\thepage\hfill
                                        {\it R. Berlanga and Ma.\hspace{1pt}\'A. Sandoval-Romero}}}
\def\@oddhead{\hbox to \textwidth{\footnotesize{\it Perturbations and Deformations}\hfill\thepage}}
\catcode`@=12

\def\qed{\unskip \enskip \null \nobreak \hfill {\bf q.e.d.} \par}

\begin{document}
 
\thispagestyle{empty}

\begin{flushleft}
    {\large\bf Small Perturbations and Infinitesimal Deformations 

on Surfaces of Revolution}
\end{flushleft}
  
\vspace*{2 true mm}
\begin{flushleft}
{\bf Ricardo Berlanga\\
\vspace*{2 true mm}
                {\small\it Instituto de Investigaciones en Matem\'aticas Aplicadas y en Sistemas (IIMAS) \\
                Departamento de F\'isica Matem\'atica \\   
                Universidad Nacional Aut\'onoma de M\'exico (UNAM) \\
                04510 M\'exico D.F. \  M\'exico\\
                E-mail: berlanga@servidor.unam.mx}}
\end{flushleft}                

\begin{flushleft}
{\bf Mar\'ia de los \'Angeles Sandoval Romero\\
\vspace*{2 true mm}
                {\small\it Facultad de Ciencias\\
                Departamento de Matem\'aticas \\
                Universidad Nacional Aut\'onoma de M\'exico (UNAM) \\
                04510 M\'exico D.F. \  M\'exico\\
                E-mail: selegna@ciencias.unam.mx}}
\end{flushleft}                
 
\vspace*{3 true mm}
\noindent{\small {\small\bf Abstract.} \ \ 
Any given surface of revolution embedded in Euclidean three-space can always be perturbed by arbitrarily small ambient isotopies 
as to admit highly nontrivial vector fields inducing infinitesimal deformations. 
For this matter Morse Theory is used, clarifying and giving a generalization of a problem originaly introduced by M. Spivak \cite{Spi01} in a modern perspective.}
\vspace{3mm}

\noindent{\small\it 2010 Subject Classification.}\ \  53A05, 53C05\\
 

\section{Introduction}
\vspace*{2 true mm}

There are many different ways in which a geometric structure may be rigid, unwarpable, bendable or deformable. A precise formulation of 
any of these intuitions may be, for example, macroscopic or infinitesimal and may involve the notion of an (ambient) isotopy
or not. The interest in this type of problem matured in the late part of the nineteenth century (see \cite{Str01}) and
since then, the vast accumulated literature on the subject has also covered all sorts of applications in Architecture, Engineering and 
Mechanics (see  \cite{Sab01}, \cite{Sab02}, \cite{Sab03}, \cite{Pog01}). In any event, the paper of Stefan Cohn-Vossen, published in 1930, 
(see \cite{Coh01} ) is certainly a classic approach to our main subject.

From the mathematical point of view, the partial differential equations and the analytic methods involved have been a matter of recent 
interest (see \cite{Yau01} and \cite{Yau02}). Despite the wealth of important results in the area, Shing-Tung Yau states, for example,  
that ``the study of the rigidity of nonconvex surfaces is still in its infancy'' (see \cite[p. 27]{Yau01}).  

The main result of this work is:

\medskip

 \noindent  {\sl Let $\lambda = (r, h) : S^1 \rightarrow {\mathbb R}^2$ be an embedding such 
that $r(s) > 0, \ \forall s \in S^1$, and let ${\mathcal K}  = \lambda(S^1)$. 
Then, for any compact set $B \subseteq {\mathbb R}^2$ containing ${\mathcal K}$ in its interior, there is an ``admissible''
perturbation $\alpha : [o,1] \times {\mathbb R}^2 \rightarrow {\mathbb R}^2$ supported in $B$ such that 
the surface of revolution obtained by revolving the perturbation of $\lambda$ at time one by $\alpha$ admits a nontrivial infinitesimal deformation.}

\medskip

This theorem starts with an idea of M. Spivak (see \cite[Vol. 5, pp. 253-260]{Spi01}) about 
the infinitesimal deformations of surfaces of revolution. Spivak works with a particular kind of surface 
diffeomorphic to the two dimensional sphere satisfying a restrictive analytical condition at its poles. 
In contrast, our methods work for general surfaces of revolution, 
diffeomorphic either to the 2-sphere or to the 2-dimensional torus.
For this matter, the elementary part of Morse Theory is used, thus introducing a new way of combining analytical and geometric methods 
which allows substantial generalization and clarification. This is done in section 3, where further remarks are given. 

\label{intro}


\section{Infinitesimal Rigidity}
\noindent {\bf Definition.} {\sl Let $N$ be a manifold and let $\varphi : N \rightarrow {\mathbb R}^m$ be an embedding. A {\bf bending} 
of $\varphi$ is a smooth map $\alpha : [ 0, 1 ] \times N \rightarrow {\mathbb R}^m$ such that
\vspace*{2 true mm}

(a) each map $\bar{\alpha} (t) : N \rightarrow {\mathbb R}^m$, given by $p \mapsto \alpha (t, p)$, is an embedding,

(b) $\bar{\alpha} (0) = \varphi$

(c) for all $t \in [ 0, 1 ]$, the pull-back $\bar{\alpha} (t)^*( < , > )$ equals $\varphi^*( < , > )$.

\vspace*{2 true mm}
 
The bending $\alpha : [ 0, 1 ] \times N \rightarrow {\mathbb {R}}^m$ is called {\bf trivial} if each  
$\bar{\alpha} (t)$ is $A_t \circ \varphi$ for some isometry $A_t$ of the ambient space ${\mathbb R}^m$; 
it is called {\bf nontrivial} if at least one  $\bar{\alpha} (t)$ is not of this form. We say that the
embedding $\varphi$ is {\bf bendable} if there is a nontrivial bending of $\varphi$; otherwise it is called 
{\bf unbendable}.}

\vspace*{2 true mm}

$N$ be a manifold and let $TN$ denote its tangent bundle. For $G : N \rightarrow {\mathbb{R}}^m$, 
let $TG : TN \rightarrow T\mathbb{R}^m \cong {\mathbb{R}}^m \times {\mathbb{R}}^m$ be its derivative and
let  $\mathrm{d}G : TN \rightarrow {\mathbb{R}}^m$ be the projection of $TG$ on the second coordinate. 
Let $U \subseteq {\mathbb{R}}^k$ be open an let $G : U \rightarrow {\mathbb{R}}^m$ be smooth. 
Denote by $\mathrm{J} G : U \times {\mathbb{R}}^k \rightarrow {\mathbb{R}}^m$ the Jacobian matrix of $G$. 
Endow ${\mathbb{R}}^m$ with the Riemannian metric defined by its dot product $\bullet$.
If $A$ is a matrix, denote its transpose by $A^t$.

\vspace*{2 true mm}

\noindent {\bf Definition.} {\sl Let $\varphi : N \rightarrow {\mathbb R}^m$ be an embedding. A vector field 
$Z : N \rightarrow {\mathbb R}^m$ is an infinitesimal deformation for $\varphi$ if it satisfies: 
}
$$
\mathrm{d}Z(X) \bullet \mathrm{d}\varphi(X) \ =\ 0 \ \ \ \ \ \forall\ X \in \mathfrak X( N )
$$

\vspace*{2 true mm}

Infinitesimal deformations arising from trivial bendings are called {\bf trivial}.
If $\varphi$ admits only trivial infinitesimal deformations, then it is called {\bf infinitesimally rigid}.

\vspace*{2 true mm}

It is not difficult to verify the following proposition.

\vspace*{2 true mm}

\noindent {\bf  Proposition.} {\sl  Let $N$ be a manifold of dimension $k$ and let 
$\varphi : N \rightarrow {\mathbb R}^m$ be an embedding. A vector field $Z : N \rightarrow {\mathbb R}^m$ 
is an infinitesimal deformation for $\varphi$ if and only if
there is an atlas $\mathcal{A}\ =\ \{ (U, \beta) \mid \beta : U \rightarrow {\mathbb R}^k \}$ of $N$ such that, 
for each $\beta \in \mathcal{A}$, the matrix 
$[\mathrm{J}(Z \circ \beta^{-1})]^t \cdot \mathrm{J}(\varphi\circ \beta^{-1})$ is antisymmetric.}

\vspace*{2 true mm}

\label{sec:2}


\section{Perturbations of Surfaces of Revolution}

Let $S^1$ be the unit circle thought of as the real numbers modulo $2\pi$ 
(ie. $s \sim ( \cos s, \sin s ), \forall s \in \mathbb R)$. Therefore we may think of a circle point as a unit complex number or
as a lateral class $s + 2\pi \cdot {\mathbb Z}$. 
Let $\lambda = (r, h) : S^1 \rightarrow {\mathbb R}^2$ be an embedding such that $r(s) > 0, \ \forall s \in S^1$.
Set $\gamma (t)  = (\cos t, \sin t, 0)$ and $\widehat k = (0, 0, 1)$.
Let $F : S^1 \times {\mathbb R} \rightarrow {\mathbb R}^3$ be given by 
$F(s,t) \ =\ r(s) \cdot \gamma(t)\ +\ h(s) \cdot \widehat k$. 
Therefore,  $\mathcal {S} = F(S^1 \times {\mathbb R})$ is a smooth surface and $F$ is a local embedding. 
With some ambiguity we make the following definition.

\vspace*{2 true mm}

\noindent {\bf Definition.} {\sl $\mathcal {F}(\lambda) = F$ (or $\mathcal {S}$) is called the {\bf surface 
of revolution obtained by revolving the curve $\lambda$ around the vertical axis and (globally) parametrized 
by $F$}}.

\vspace*{2 true mm}  

\noindent {\bf Definitions.} {\sl Let $N$ be a topological space and set $l = 0$; 
or let $N$ be a smooth manifold and set $l = 3$. 
A  $C^l$ map $\alpha : [ 0, 1 ] \times N \rightarrow N $ is an {\bf isotopy} of $N$ if
(a) each map $\bar{\alpha}(t): N \rightarrow N$, given by $p \mapsto \alpha (t, p)$, 
is a $C^l$ homeomorphism, (b) $\bar{\alpha} (0) = Id : N \rightarrow N$.}

\vspace*{2 true mm}

{\sl Let $\alpha : [ 0, 1 ] \times {\mathbb R}^2 \rightarrow {\mathbb R}^2$ 
be an isotopy of ${\mathbb R}^2$ and let ${\mathcal K} \subseteq \{ (x, z) \in  {\mathbb R}^2 \mid x > 0 \}$. 
We say that $\alpha$ is an {\bf admissible perturbation for} ${\mathcal K}$
if $\bar{\alpha} (t) ({\mathcal K}) \subseteq \{ (x, z) \in  {\mathbb R}^2 \mid x > 0 \}$ for all $t \in [ 0, 1 ]$}.

\vspace*{2 true mm}

{\sl Let $G : N \rightarrow N$ be a $C^l$ map. The closure of  the set 
$\{ p \in N \mid G(p) \neq  p \}$ is the {\bf support} of $G$ and is denoted by supp($G$). 
For an isotopy $\alpha$ define its {\bf support} as the closure of the set 
$\bigcup\ \{ \mathrm{supp} \  \bar{\alpha} (t) \mid t \in [ 0, 1 ] \}$.}

\vspace*{2 true mm}

The next statement is the main result of this work.

\vspace*{2 true mm}

\noindent {\bf Theorem.} {\sl Let $\lambda = (r, h) : S^1 \rightarrow {\mathbb R}^2$ be an embedding such 
that $r(s) > 0, \ \forall s \in S^1$, and let ${\mathcal K}  = \lambda(S^1)$. 
Then, for any compact set $B \subseteq {\mathbb R}^2$ containing ${\mathcal K}$ in its interior, there is an admissible
perturbation $\alpha$ for ${\mathcal K}$ supported in $B$ such that $\mathcal {F}(\bar{\alpha} (1) \circ \lambda)$
admits a nontrivial infinitesimal deformation.}

\vspace*{2 true mm}

\noindent {\bf Definition.} {\sl Let $N$ be a manifold and let $g : N \rightarrow {\mathbb R}$ be smooth. 
A point $p \in N$ is {\bf critical} if $\mathrm{d}g\mid_p\ = 0$. Let ${\mathcal C}(g)$ be the set of critical points of 
$g$ in $N$. A critical point $p$ is {\bf nondegenerate} if, for some (hence any) coordinate chart $\beta$ around $p$, 
the Hessian matrix of $g \circ \beta^{-1}$ at $\beta(p)$ is nonsingular. 
Functions whose critical points are nondegenerate are called {\bf Morse functions}. }   

\vspace*{2 true mm}

For an arbitrary oriented one dimensional manifold ${\mathcal K} \subseteq  {\mathbb R}^2$, 
let $\kappa : {\mathcal K} \rightarrow {\mathbb R}$  denote its (signed) curvature and 
let $\zeta : {\mathcal K} \rightarrow  {\mathbb R}$ be the restriction to 
${\mathcal K}$ of the projection $(x,z) \mapsto z$.

\vspace*{2 true mm}

\noindent {\bf Lemma 1.} {\sl Let ${\mathcal K} \subseteq  {\mathbb R}^2$ be a compact one dimensional manifold.
Suppose the height function $\zeta : {\mathcal K} \rightarrow {\mathbb R}$ is Morse.
Then

(a) ${\mathcal C}(\zeta)$ is finite,

(b) $\kappa(p) \not = 0$ for all $p \in {\mathcal C}(\zeta)$.    

(c) For each $p \in {\mathcal C}(\zeta)$, there is a coordinate chart $\beta : U \rightarrow\ <-\epsilon, \epsilon>$ 
around $p = \beta^{-1}(0)$ such that $\zeta \circ \beta^{-1}(u) = \zeta (p) \pm u^2$, 
for all $u \in\ <-\epsilon, \epsilon>$.}
 
\noindent {\bf {\textsc Proof.}} 
For any Morse function the set of critical points is closed and discrete (see Hirsch \cite[p. 143]{Hir01}), 
hence ${\mathcal C}(\zeta)$ is finite, for ${\mathcal K}$ is compact. 

Let $p \in {\mathcal C}(\zeta)$. By Morse's Lemma (\cite[p. 145]{Hir01}), 
there is a coordinate chart $\beta : U \rightarrow\ <-\epsilon, \epsilon>$ 
around $p = \beta^{-1}(0)$ such that $z(u) = \zeta (p) + \sigma \cdot u^2$, for all $u \in\ <-\epsilon, \epsilon>$,
where $z = \zeta \circ \beta^{-1}$ and $\sigma = 1$ if the index of $p$ is zero and $\sigma = -1$ 
if the index of $p$ is one. 
Let $x :  <-\epsilon, \epsilon> \rightarrow {\mathbb R}$ be the first coordinate of 
$\beta^{-1}$, so $\beta^{-1} = (x, z)$. Since $z^\prime (0) = 0$,
and $\beta$ is a chart, $x^\prime (0) \not = 0$. Let $\phi :  <-\epsilon, \epsilon> \rightarrow <-\epsilon, \epsilon>$
be such that $\phi (u) = -u$.
If $\beta$ is nonoriented, substitute $\beta$ for $\beta \circ \phi$, 
so we may always assume that $\beta$ is oriented if necessary. In any case, $z$ is unchanged.
Now, 
$$
\kappa(p) =  \frac{\mathrm{det}((\beta^{-1})^{\prime}(0), (\beta^{-1})^{\prime \prime}(0))}
{\| (\beta^{-1})^{\prime}(0) \|^3}
=  \frac{x^\prime (0)\cdot z^{\prime \prime}(0)}{\| (\beta^{-1})^{\prime}(0) \|^3}
=  \frac{x^\prime (0)\cdot 2 \sigma}{\| (\beta^{-1})^{\prime}(0) \|^3}
$$ 
\noindent showing that $\kappa(p)\not = 0$ \qed

\vspace*{2 true mm}

\noindent {\bf Assertion.} {\sl Let ${\mathcal K} \subseteq  {\mathbb R}^2$ be a compact one dimensional manifold.
For each $\nu \in S^1$ let $f_\nu : {\mathcal K} \rightarrow {\mathbb R}$ be the map $f_\nu(p) = \nu \bullet p$,
so the orthogonal projection into the line through $\nu$ is given by $p \mapsto f_\nu(p) \cdot \nu$. 
Then the set of $\nu \in S^1$ such that $f_\nu$ is a Morse function is open and dense.}

\noindent {\bf {\textsc Proof.}} This is a minor instance of Hirsch \cite[Ex. 2, p. 148]{Hir01} \qed

\vspace*{2 true mm}

\noindent {\bf Corollary.} {\sl Let ${\mathcal K} \subseteq  {\mathbb R}^2$ be a compact one dimensional manifold.
Let $R_\theta$ denote the plane rotation at angle $\theta$ in the positive direction with center at the origin
and let ${\mathcal K}_\theta = R_\theta(\mathcal K)$.
Then, for almost every $\theta \in {\mathbb R}$, the height function 
$\zeta_\theta : {\mathcal K_\theta} \rightarrow {\mathbb R}$ is Morse.}

\noindent {\bf {\textsc Proof.}}   Let $\widehat \iota = (0,1)$; the correspondence 
$\theta \mapsto R_{- \theta}(\ \widehat \iota\ )$ is a 
continuous bijection from $[0, 2\pi>$ onto $S^1$; let $\widehat \iota = R_\theta(\nu)$. 
Then, $f_\nu  = \zeta_\theta \circ R_\theta$
for $\zeta_\theta (R_\theta(p)) =$ $\widehat \iota \bullet R_\theta(p) =$ $R_\theta(\nu) \bullet R_\theta(p) =$ 
$\nu \bullet p = f_\nu (p)$. 
Since $R_\theta : {\mathcal K} \rightarrow {\mathcal K}_\theta$ is a diffeomorphism, then $f_\nu$ is
Morse iff  $\zeta_\theta$ is \qed

\vspace*{2 true mm}

\noindent {\bf Definition.}  {\sl Let ${\mathcal K} \subseteq  {\mathbb R}^2$ be a one dimensional manifold with Morse
height function $\zeta$ and let $p \in {\mathcal C}(\zeta)$. Let $x$ be the horizontal coordinate $(x,z) \mapsto x$.
Then $x$ must be a diffeomorphism around $p$; therefore, $\zeta \circ x^{-1}$ is locally defined. 
The point $p$ is called {\bf even analytic} if  $\zeta \circ x^{-1}$ is a real even analytic function around $x(p)$. 
${\mathcal K}$ is {\bf even analytic (at ${\mathcal C}(\zeta)$)} if all its (nondegenerate) critical points 
are even analytic.} 

\vspace*{2 true mm}

\noindent {\bf Lemma 2.} {\sl ${\mathcal K} \subseteq \{ (x,z) \in {\mathbb R}^2 \mid x > 0 \}$ 
be a compact one dimensional manifold.
Then, for any compact set $B \subseteq {\mathbb R}^2$ containing ${\mathcal K}$ in its interior, there is an admissible
perturbation $\alpha$ for ${\mathcal K}$ supported in $B$ and arbitrarily close to the identity such that
$\bar{\alpha} (1) ({\mathcal K})$ has Morse height function $\zeta_1$ and is even analytic 
(at ${\mathcal C}(\zeta_1)$)}.

\noindent {\bf {\textsc Proof.}}   
Let $B \subseteq {\mathbb R}^2$ be compact containing ${\mathcal K}$ in its interior $B^\circ$. 
Assume, without loss of generality that $B \subseteq \{ (x,z) \in {\mathbb R}^2 \mid x > 0 \}$.
Let $C \subseteq {\mathbb R}^2$ be such that ${\mathcal K} \subseteq C^\circ \subseteq C \subseteq B^\circ$.
Let $\eta : {\mathbb R}^2 \rightarrow [0, 1]$ be a smooth (bump) function with $\eta \mid_C\ = 1$ and   
$\eta \mid_{{\mathbb R}^2\setminus B^\circ}\ = 0$. Let $X$ be the vector plane field given by 
$X(x,z) = \eta(x,z) \cdot (-z, x), \forall (x,z) \in {\mathbb R}^2$. The one parameter group of diffeomorphisms 
$R : {\mathbb R} \times {\mathbb R}^2 \rightarrow {\mathbb R}^2$ induced by $X$ is such that every $\bar R(t)$ is
supported in $B$ and there exists $\epsilon > 0$ where $R(t,p)$ is the rotation at angle $t$ in the positive
direction with center at the origin ($\forall\ \mid t \mid <\epsilon$ and $\forall\ p \in {\mathcal K}$).  

Choose $0 < \theta < \epsilon$ such that the height function 
$\zeta_\theta : \bar R(\theta)({\mathcal K}) \rightarrow {\mathbb R}$ is Morse. Note that $\theta$ may be 
chosen as small as desired ($\epsilon$ is not necessarily small). 

Let $p = (x_0, z_0) \in {\mathcal C}(\zeta_\theta)$ be a concave critical point. So choose $\delta_1  > 0$ 
such that

\noindent (1) $\rho = \zeta_\theta \circ x^{-1}$ is defined on $<x_0 - 2 \delta_1, x_0 + 2 \delta_1>$, 

\noindent (2) $\rho^\prime$ is strictly decreasing on $<x_0 - \delta_1, x_0 + \delta_1>$.

\vspace*{2 true mm}
By making $\delta_1$ smaller, elementary Taylor theory says that there exist $K < 0$ such that 
\vspace*{2 true mm}

\noindent (3) $\rho(x) \leqslant z_0 + K \cdot (x - x_0)^2$ for all $x \in <x_0 - \delta_1, x_0 + \delta_1>$.

\vspace*{2 true mm}
By making $\delta_1$ smaller, now we can find $\delta_2 > 0$ such that
\vspace*{2 true mm}

\noindent (4) $[x_0 - \delta_1, x_0 + \delta_1] \times [z_0 - \delta_2, z_0 + 3 \delta_2 / 2] \subseteq B^\circ$,  

\noindent (5) $z_0 - \delta_2 / 2 < \rho (x)$ for all $x \in <x_0 - \delta_1, x_0 + \delta_1>$,

\noindent (6)  $(<x_0 - \delta_1, x_0 + \delta_1> \times <z_0 - \delta_2, z_0 + 3 \delta_2 / 2> ) \cap \bar R(\theta)({\mathcal K})$
is equal to $\{ (x,\rho(x)) \mid x \in <x_0 - \delta_1, x_0 + \delta_1>  \}$.

\vspace*{2 true mm}
Let $\eta_1 : {\mathbb R} \rightarrow [0,1]$ be a smooth function such that
\vspace*{2 true mm}

\noindent (7) $\eta_1 = 1$ on  $[x_0 - \delta_1 / 2, x_0 + \delta_1 / 2]$,

\noindent (8) $\eta_1 = 0$ out of  $<x_0 - \delta_1 , x_0 + \delta_1 >$,

\noindent (9) $\eta_1$ is increasing on  $<x_0 - \delta_1 , x_0 - \delta_1 / 2 >$ and decreasing on
                                               $<x_0 + \delta_1 / 2, x_0 - \delta_1 >$. 

\vspace*{2 true mm}
Let $\eta_2 : {\mathbb R} \rightarrow [0,1]$ be a smooth function such that
\vspace*{2 true mm}

\noindent (10) $\eta_2 = 1$ on  $[z_0 - \delta_2 / 2, z_0 + \delta_2 ]$,

\noindent (11) $\eta_2 = 0$ out of  $<z_0 - \delta_2 , z_0 + 3 \delta_2 / 2 >$.

\vspace*{2 true mm}

Define $f(x) = z_0 + K (x - x_0)^2$, which is certainly even analytic at $x_0$, and let $Y$ be the vector 
field given by 
$Y(x,z) = (0, \eta_2(z) \eta_1(x)  (f(x) - \rho(x)), \forall (x,z) \in {\mathbb R}^2$, (let $Y(x,z) = (0,0)$ 
when $\rho(x)$ fails to exist). 
Note that $Y$ is vertical and points upwards when is different than zero.
The one parameter group of diffeomorphisms 
$S : {\mathbb R} \times {\mathbb R}^2 \rightarrow {\mathbb R}^2$ induced by $Y$ is such that every $\bar S(t)$ is
supported in $[x_0 - \delta_1, x_0 + \delta_1] \times [z_0 - \delta_2, z_0 + 3 \delta_2 / 2]$. 
The flow for $Y$ in the box $[0,1] \times {\mathbb R} \times <z_0 - \delta_2 / 2 , z_0 + \delta_2 >$ is given by
$(t,(x, z)) \mapsto (x, z + t \eta_1(x) (f(x) - \rho(x)))$.
For $(t,x) \in [0,1] \times <x_0 - \delta_1, x_0 + \delta_1>$, (and condition (5)), we have 
$(t,(x, \rho(x))) \mapsto (x, \rho(x) + t \eta_1(x) (f(x) - \rho(x)) = (x, g_t(x))$.
In time one, for $x \in <x_0 - \delta_1 / 2, x_0 + \delta_1 / 2 >$, we have $g_1 (x) = f(x)$.  
Now, $g_t^\prime (x) = (1 - \eta_1(x))\ \rho^\prime(x) + \eta_1(x)\ f^\prime(x) + \eta^\prime_1(x)\ (f(x) - \rho(x))$.
By (2),(3),(7),(8) and (9) above  $g_t^\prime (x)$ has $x_0$ as its only zero. 
Therefore, by (6), the perturbation $ \bar R(\theta)({\mathcal K}) \rightarrow \bar S(t)(\bar R(\theta)({\mathcal K}))$ 
does not introduce new critical points and $\bar S(1)(\bar R(\theta)({\mathcal K}))$ is even analytic at $p$.

The situation for convex critical points is treated in exactly the same way.  Hence, we can perturb 
$\bar R(\theta)({\mathcal K})$
at all points in ${\mathcal C}(\zeta_\theta)$ (disjointly) in order to make all resulting critical points even analytic
\qed

\vspace*{2 true mm}

Lemma 2 and the following Proposition imply the Theorem of this section.

\vspace*{2 true mm}

\noindent {\bf Proposition.} {\sl Let ${\mathcal K} \subseteq \{ (x,z) \in {\mathbb R}^2 \mid x > 0 \}$ be diffeomorphic 
to the unit circle $S^1$.
Suppose ${\mathcal K}$ has Morse height function $\zeta$ which is also even analytic (at ${\mathcal C}(\zeta)$).
Let $\lambda = (r, h) : S^1 \rightarrow {\mathcal K}$ be a diffeomorphism, 
so the surface of revolution $\mathcal{F}(\lambda)$ is defined.
Then, for any compact set $B \subseteq {\mathbb R}^2$ containing ${\mathcal K}$ in its interior, there is an admissible
perturbation $\alpha$ for ${\mathcal K}$ supported in $B$ such that $\mathcal{F}(\bar{\alpha} (1) \circ \lambda)$
admits nontrivial infinitesimal deformations different from zero almost everywhere}.

\noindent {\bf {\textsc Proof.}} (cf Spivak \cite[Vol. 5, pp. 253-260]{Spi01}).  
Let $F = \mathcal {F}(\lambda)$. Then, any vector field $Z$ along $F$ can be written uniquely as 
$$
Z(s, t) \ =\ a(s, t) \cdot \gamma(t) \ +\ b(s, t) \cdot \gamma^{\ \prime}(t) \ +\ c(s, t) \cdot \widehat k
$$
for some smooth functions $a, b, c$. An infinitesimal deformation $Z$ for $F$ is trivial if it is of the form 
$Z(s,t) = \mu \times F(s,t) + \omega$, 
where $\mu, \omega \in {\mathbb R}^3$ and $\times$ denotes the cross product.  

Since $F$ is a global regular parametrization, its local inverses form an atlas for its image.
Hence, by the Proposition of section 2, $Z$ is an infinitesimal deformation if an only if
$(J Z)^t \cdot JF$ is a $2 \times 2$ antisymmetric matrix. That is,
$$
\frac{\partial Z}{\partial s} \bullet \frac{\partial F}{\partial s} = 0,\ \ \ \
\frac{\partial Z}{\partial t} \bullet \frac{\partial F}{\partial t} = 0,\ \ \ \
\frac{\partial Z}{\partial t} \bullet \frac{\partial F}{\partial s}   +
\frac{\partial Z}{\partial s} \bullet \frac{\partial F}{\partial t} = 0\ .
$$

\vspace*{2 true mm}
 
A simple calculation shows that these equations are explicitly given by

$$
r^\prime \displaystyle\frac{\partial a}{\partial s} + h^\prime \displaystyle\frac{\partial c}{\partial s} = 0,\ \ \
\displaystyle\frac{\partial b}{\partial t} + a = 0,\ \ \
r^\prime \Bigg ( \displaystyle\frac{\partial a}{\partial t} - b \Bigg  ) +
r \displaystyle\frac{\partial b}{\partial s} +
h^\prime \displaystyle\frac{\partial c}{\partial t} = 0\ \ \ (\dagger)
$$

Propose standar solutions of the form
$$
a(s,t) = \sum_{k=-\infty}^\infty e^{ikt}\phi_k(s), b(s,t) = \sum_{k=-\infty}^\infty e^{ikt}\psi_k(s), 
c(s,t) = \sum_{k=-\infty}^\infty e^{ikt}\xi_k(s)
$$

For a (complex-valued) solution involving a single $k$, equations ($\dagger$) become  

$$
\left .
\begin{array}{r c l}
r^\prime \phi_k^\prime + h^\prime \xi_k^\prime &=& 0 \\[7pt]
i k\ \psi_k + \phi_k &=& 0\\[7pt]
r^\prime ( i k\ \phi_k - \psi_k) + r \psi_k^\prime + i k\ h^\prime \xi_k &=& 0
\end{array} \right \}\ \ \ \ \ (\dagger \dagger)
$$

\vspace*{2 true mm}

If this system is solved, then it is easy to see that $a_k(s,t) = e^{ikt}\phi_k(s) + e^{-ikt}\bar\phi_k(s)$, 
$b_k(s,t) =  e^{ikt}\psi_k(s) + e^{-ikt}\bar \psi_k(s)$, $c_k(s,t) =  e^{ikt}\xi_k(s) + e^{-ikt}\bar \xi_k(s)$
are the scalar components of a real infinitesimal deformation.

Differentiating the third equation in $(\dagger \dagger)$ and using the other two, we obtain the equivalent system 

$$
\left .
\begin{array}{r c l}
\phi_k + i k\ \psi_k &=& 0 \\[7pt]
ik\ h^\prime \xi_k + (k^2 - 1) r^\prime \psi_k + r \psi_k^\prime  &=& 0\\[7pt]
r \psi_k^{\prime\prime} + (k^2 - 1) r^{\prime\prime} \psi_k + ik\ h^{\prime\prime}\xi_k &=& 0
\end{array} \right \}\ \ \ \ \ (\dagger \dagger \dagger)
$$

\vspace*{2 true mm}

System $(\dagger \dagger \dagger)$ is easily transformed into system

$$
\left .
\begin{array}{r c l}
\phi_k + i k\ \psi_k &=& 0 \\[7pt]
ik\ h^\prime \xi_k + (k^2 - 1) r^\prime \psi_k + r \psi_k^\prime  &=& 0\\[7pt]
r \Bigg [ \displaystyle\frac{\psi_k^\prime}{h^\prime} \Bigg ]^\prime + 
                   (k^2 - 1) \Bigg [ \displaystyle\frac{r^\prime}{h^\prime} \Bigg ]^\prime \psi_k &=& 0
\end{array} \right \}\ \ \ \ \ (\dagger \dagger \dagger \dagger)
$$

\vspace*{2 true mm}

Solutions of the third equation above immediately give solutions of the first two. Hence we can define infinitesimal
deformations away from the critical set $\{ F(s,t) \mid \  h^\prime(s) = 0 \}$. 
Now, $\zeta$ is given by $(x,z) \mapsto z$ and if we let $x$ to be the horizontal coordinate $(x,z) \mapsto x$, 
then $r = x \circ \lambda$ and $h = \zeta \circ \lambda$. 
Also, ${\mathbb R}$ is (periodically) partitioned into 
\vspace*{2 true mm}

(a) open intervals where $h^\prime$ is nonzero and,

\vspace*{2 true mm}

(b) its endpoints, which correspond to the elements in ${\mathcal C}(\zeta) \subseteq {\mathbb R}^2$.

\vspace*{2 true mm}
We shall refer to the closures of these intervals (or its images under $h$ or $\lambda$) as ``the segments''. 
In particular, for any $s \in {\mathbb R}$, if $h^\prime(s)=0$ then $r$ is locally invertible around $s$. 
For such points let $\rho = h \circ r^{-1} = \zeta \circ x^{-1}$, so  $\rho$ is even analytic. 
Incidentally, this proves that the number of elements in ${\mathcal C}(\zeta)$, say  
$\sharp \mid {\mathcal C}(\zeta) \mid$, is even, because intervals where $h^\prime$ is positive 
and negative must alternate, for $\lambda^{-1}({\mathcal C}(\zeta))$ contains no inflection points. 
Clearly $2 \leq \sharp \mid {\mathcal C}(\zeta) \mid$, for the top and bottom points in ${\mathcal K}$ are critical.
Also, the equality $h \circ r^{-1} = \zeta \circ x^{-1}$ shows that if 
$s_\flat  = s_\natural  + 2\cdot \pi$ correspond to the same
point in ${\mathcal C}(\zeta)$ then the associated $\rho$ is the same.
Observe that 
$\{ s \in {\mathbb R} \mid r^\prime(s) = 0\}$ may be rather nasty. 

Define, when meaningful,
$\widetilde \Phi_k = \phi_k \circ h^{-1}$, $\widetilde \Psi_k = \psi_k \circ h^{-1}$, 
                                           $\widetilde \Xi_k = \xi_k \circ h^{-1}$ 
and 
           $\Phi_k = \phi_k \circ r^{-1}$,             $\Psi_k = \psi_k \circ r^{-1}$,
                                                       $\Xi_k = \xi_k \circ r^{-1}$,
which in turn satisfy the respective systems

$$
\left .
\begin{array}{r c l}
\widetilde \Phi_k + i k\ \widetilde \Psi_k &=& 0 \\[7pt]
ik\ \widetilde \Xi_k +   
(r \circ h^{-1}) \widetilde \Psi_k^\prime  
+ (k^2 - 1) (r \circ h^{-1})^\prime \widetilde \Psi_k &=& 0\\[7pt]
(r \circ h^{-1}) \widetilde \Psi^{\prime\prime}_k + (k^2 - 1) (r \circ h^{-1})^{\prime\prime} \widetilde \Psi_k &=& 0
\end{array} \right \}\ \ \ \ \ (5-\dagger)
$$

\vspace*{2 true mm}

\noindent and \ \ $\Phi_k + i k\ \Psi_k = 0\ $,  
$ik\ \rho^\prime\ \Xi_k + \textrm{Id}\  \Psi_k^\prime  + (k^2 - 1) \Psi_k = 0\ $, 
$\textrm{Id}\ \rho^\prime  \Psi^{\prime\prime}_k -  \textrm{Id}\ \rho^{\prime\prime}\Psi^{\prime}_k
               - (k^2 - 1) \rho^{\prime\prime}\Psi_k  = 0$.
Finally, the second equation of this last system  may be changed to obtain the equivalent system 

$$
\left .
\begin{array}{r c l}
\Phi_k + i k\ \Psi_k &=& 0 \\[7pt]
ik\ \Xi_k + 
\Bigg ( \displaystyle\frac{\rho^\prime}{\rho^{\prime\prime}} \Bigg )\textrm{Id}\ \Psi_k^{\prime\prime}
&=& 0\\[14pt]
\textrm{Id}\ \rho^\prime  \Psi^{\prime\prime}_k -  \textrm{Id}\ \rho^{\prime\prime}\Psi^{\prime}_k
               - (k^2 - 1) \rho^{\prime\prime}\Psi_k  &=& 0
\end{array} \right \}\ \ \ \ \ (6-\dagger)
$$

Note that, when both $r^{-1}$ and $h^{-1}$ are defined, then $\widetilde \Phi_k = \Phi_k \circ \rho^{-1}$, 
$\widetilde \Psi_k = \Psi_k \circ \rho^{-1}$ and $\widetilde \Xi_k = \Xi_k \circ \rho^{-1}$. Again, solutions 
for the third equation in $(5-\dagger$) immediately produce solutions for the first two away from critical points.
In the same manner, solutions for the third equation in $(6-\dagger$) immediately produce solutions for the first 
two around critical points.

Fix $p \in {\mathcal C}(\zeta)$. Then, $\rho$ is (even) analytic around $x_0 = x(p)$ with zero derivative at $x_0$.
Then it can be expressed as $\rho (x) = \zeta(p) + (1/2!) \rho^{\prime\prime}(x_0) (x - x_0)^2 +
(1/4!)\rho^{\prime\prime\prime\prime}(x - x_0)^4 + \cdots$,
with  $\rho^{\prime\prime} (x_0) \not = 0$. Hence we can write $\rho\prime (x) =  (x - x_0) P(x)$ with 
$P(x_0) = \rho^{\prime\prime} (x_0)$.  Multiplying the third equation in $(6-\dagger$) by
$(x - x_0) / x P(x)$ we get

{\footnotesize
$$
(x - x_0)^2 \Psi^{\prime\prime}_k (x) + (x - x_0) \Bigg ( - \displaystyle\frac{\rho^{\prime\prime}(x)}{P(x)} \Bigg )
\Psi^\prime_k(x) + \Bigg (-(k^2-1)\displaystyle\frac{(x-x_0)\rho^{\prime\prime}(x)}{x P(x)} \Bigg )\Psi_k (x) = 0
$$
}

This is a standard expression of a second-order differential equation having a regular singular point at $x_0$
with roots $0$ and $2$ for its ``indicial equation''. Therefore, it has a unique analytic solution of the form 
$\Psi_k (x) = (x - x_0)^2 A(x)$, satisfying $A(x_0) = 1$, (see Birkhoff and Rota \cite[p. 252 ff]{Bir01}). 
In particular, $x_0$ is an isolated zero for $\Psi_k (x)$. 
Note that $x_0 >0$, but this is a geometric condition that assures that we have a smooth surface of revolution.
From the analytic point of view, $x_0 = 0$ is admisible. In this case the indicial equation has roots $1 \pm k$.
If $k > 1$, there is  exactly one analytic solution (around $0$) of the form $\Psi_k (x) = x^{1 + k} A(x)$, 
(cf \cite[Vol. 5, p. 257]{Spi01}).  

Let $s_1 < s_2 < s_3$ be three consecutive points in $\lambda^{-1}({\mathcal C}(\zeta))$, where 
$s_3 = s_1 + 2 \cdot \pi$ is possible. Then there is an $\epsilon > 0$ such that
the system $(\dagger \dagger \dagger \dagger )$ has well determined nontrivial smooth solutions $\phi_k$, $\psi_k$,
$\xi_k$ in the interval $<s_1 - \epsilon, s_2>$. In particular, $\psi_k(s_1) = \psi^\prime_k(s_1) = 0$ and 
$\xi_k(s_1) = 0$.

It is true that if $\psi_k(s_2)$ is defined and equal to zero, then the solutions $\phi_k$, $\psi_k$,
$\xi_k$ are all well defined in the interval $<s_1 - \epsilon, s_2+ \epsilon>$ and therefore, 
in the interval $<s_1 - \epsilon, s_3>$. The obstruction to this continuation process is that, usually, we expect to
have $\psi_k(s_2) = \pm \infty$. 

Up to now $k$ has been kept fixed. To solve the above mentioned obstruction we will have to vary this integer and 
perturb ${\mathcal K}$ (a third time, but now away from ${\mathcal C}(\zeta)$).  

\vspace*{2 true mm}

So let $s_1 < s_2$ be two consecutive points in $\lambda^{-1}({\mathcal C}(\zeta))$. 
In $[s_1 , s_2]$ the function $h$ is monotone increasing or decreasing.
Let $\{ h(s_1), h(s_2)\} = \{ z_1, z_2\}$ be labeled so that $z_1 < z_2$.
Let $R = r \circ h^{-1} :\ [z_1, z_2] \rightarrow <0, \infty >$ be the horizontal coordinate for this segment.
Fix a point $z_0 \in <z_1,z_2>$ and let $x_0 = R(z_0)$. Choose $\delta_1 , \delta_2 > 0$ such that 

\noindent (1) $z_1 < z_0 - 3 \delta_2 < z_0 + 3 \delta_2 < z_2$, 

\noindent (2) $[x_0 - 2 \delta_1, x_0 + 2 \delta_1] \times [ z_0 - 3 \delta_2, z_0 + 3 \delta_2] \subseteq B^\circ$,

\noindent (3) $0 <x_0 - \delta_1/ 2 < R(z) < x_0 + \delta_1 /2,\ \forall z \in < z_0 - 3 \delta_2 ,z_0 + 3 \delta_2>$, 

\noindent (4) $[x_0 - 2 \delta_1, x_0 + 2 \delta_1] \times [ z_0 - 3 \delta_2, z_0 + 3 \delta_2] \cap {\mathcal K} 
               = \{ (R(z), z) \mid z \in [ z_0 - 3 \delta_2, z_0 + 3 \delta_2]  \}$.   
 
As in Lemma 2, we can perturb the graph of $R$ on $[ z_0 - 3 \delta_2, z_0 + 3 \delta_2]$ within the box
$[x_0 - \delta_1 / 2, x_0 + \delta_1/ 2] \times [ z_0 - 5 \delta_2 / 2, z_0 + 5 \delta_2 / 2]$, 
by a small flow $T$ induced by a horizontal vector field, in such a way that 
$[x_0 - \delta_1 / 2, x_0 + \delta_1 /2 ] \times [ z_0 - 2 \delta_2, z_0 + 2 \delta_2] \cap \bar T(1)({\mathcal K})$
is the graph of a smooth strictly concave function 
$f_0 : [ z_0 - 2 \delta_2, z_0 + 2 \delta_2] \rightarrow <x_0 - \delta_1 / 2, x_0 + \delta_1 /2 >$.  
It is clear that this perturbation does not introduce new critical points in $\bar T(1)({\mathcal K})$.
Let $\widetilde R_0 :\ [z_1, z_2] \rightarrow <0, \infty >$ denote the (new) horizontal coordinate for the segment
$[z_1, z_2]$.

Let $f_1 : [ z_0 - 2 \delta_2, z_0 + 2 \delta_2] \rightarrow <x_0 - \delta_1 / 2, x_0 + \delta_1 /2 >$ be a smooth 
function such that

\noindent (5) $f_0$ and $f_1$ have the same values and the same first and second derivatives at $z_0 \pm  2 \delta_2$,

\noindent (6) $f_1$ is strictly concave in 
   $[ z_0 - 2 \delta_2, z_0 - \delta_2> \cup < z_0 + \delta_2, z_0 + 2 \delta_2>$, 

\noindent (7) $f_1$ is strictly convex in $< z_0 - \delta_2, z_0 + \delta_2>$, 

\noindent (8) $f_1(z) \leq f_0(z),\ \  \forall z \in [ z_0 - 2 \delta_2, z_0 + 2 \delta_2]$.

The horizontal field $X(x,z) = (f_1(z) - f_0(z) , 0), \forall z \in [ z_0 - 2 \delta_2, z_0 + 2 \delta_2]$ and 
$X(x,z) = (0 , 0), \forall z \not\in [ z_0 - 2 \delta_2, z_0 + 2 \delta_2]$ is $C^2$ and may be damped by a bump 
function (in variable x) in order to produce a small flow that sends, in time one, the graph of $f_0$ onto the 
graph of $f_1$. 
Let $\{ \widetilde R_t :\ [z_1, z_2] \rightarrow <0, \infty > \}_{t \in [0,1]}$ be the continuous 
family of horizontal coordinates produced by the flow (for the segment $[z_1,z_2]$). Fix $k \in {\mathbb N}$ and 
consider the family of equations (for $t \in [0,1]$, on $<z_1, z_2>$)

$$
\widetilde \Psi^{\prime\prime}\ +\ \Bigg [\frac{(k^2-1)(\widetilde R_t)^{\prime\prime}}{\widetilde R_t}\Bigg ]\ \widetilde \Psi = 0 \ \ \ \ \ \ \ (\clubsuit) 
$$

\noindent and let $\{ \widetilde \Psi_{k, t} \}_{t \in [0,1]}$ be a continuous family of solutions 
such that its expression in the horizontal axis has $R(z_1)$ as a zero of order two. Then, 

\noindent (9) For all $t \in [0,1]$, $\widetilde R_t = \widetilde R_0$ in $[ z_1 , z_0 - 2 \delta_2]$. 
Hence, for any $t_\natural,  t_\flat  \in [0, 1]$,  $\widetilde\Psi_{k, t_\natural } \mid_{[ z_1 , z_0 - 2 \delta_2]} 
= \widetilde\Psi_{k, t_\flat } \mid_{[ z_1 , z_0 - 2 \delta_2]}$ 

\noindent (10) For all $t \in [0,1]$, $\widetilde R_t = \widetilde R_0$ in $[ z_0 + 2 \delta_2, z_2>$\ so  
$\{ \widetilde \Psi_{k, t} \}_{t \in [0,1]}$ is a family of solutions of the same differential equation in 
that interval. Therefore, by the Sturm Comparison Theorem, for any $t_\natural,  t_\flat  \in [0, 1]$,  
$-1 \leq \sharp [(\widetilde \Psi_{k, t_\natural })^{-1}(\{ 0 \})] - 
\sharp [(\widetilde \Psi_{k, t_\flat })^{-1}(\{ 0 \})] \leq 1$,
where, as before, $\sharp [V]$ denotes the cardinality of set $V$.

\noindent (11) By concavity of $f_0 = \widetilde R_0$ in $[ z_0 - 2 \delta_2, z_0 + 2 \delta_2]$,
$\widetilde \Psi_{k, 0}$ has at most one zero in the interval.

\noindent (12) By concavity of $f_1 = \widetilde R_1$ in $[ z_0 - 2 \delta_2, z_0 - \delta_2> \cup < z_0 + 
\delta_2, z_0 + 2 \delta_2]$, $\widetilde \Psi_{k, 1}$ has at most one zero in each of these intervals.

\vspace*{2 true mm} 

Considerations (9 - 12) hold for any $k$. Now, if $k$ is sufficiently large, by the Sturm Comparison Theorem  
and the convexity of $\widetilde R_1$ in $< z_0 - \delta_2, z_0 + \delta_2>$, $\widetilde \Psi_{k, 1}$ can 
be made to have many zeroes in the interval. For such value of $k$,
$\sharp [(\widetilde \Psi_{k, 0 })^{-1}(\{ 0 \})] < \sharp [(\widetilde \Psi_{k, 1 })^{-1}(\{ 0 \})]$.

Let $\hat t \in [0, 1]$ be a point of discontinuity of $ t \mapsto \sharp [(\widetilde \Psi_{k, t })^{-1}(\{ 0 \})]$.
We want to prove that $\widetilde \Psi_{k, \hat t }(z_2) = 0$. Let 
$\{ \Upsilon_{0, t}\ , \Upsilon_{\infty, t} \}_{t \in [0,1]}$ be a family of 
linearly independent solutions of equation 
$(\clubsuit )$ in $<z_1,z_2>$ such that $\Upsilon_{0, t} (z_2) = 0$ and 
$\Upsilon_{\infty, t} (z_2) \not = 0$ (possibly $\pm \infty$), for all $t \in [0,1]$. Then there exist 
unique functions $d_0 , d_\infty : [0,1] \rightarrow {\mathbb R}$ satisfying 
$\widetilde \Psi_{k, t } = d_0(t) \cdot \Upsilon_{0, t} + d_\infty(t) \cdot \Upsilon_{\infty, t}$\ .
Since the Wronskian at $z_0$ is not zero, then 
$$
\left [ \begin{array}{ c }
d_0(t)               \\ [10pt]
d_\infty(t)
\end{array} \right ] 
=  
\left [ \begin{array}{ r c }
\Upsilon_{0, t} (z_0)                  &  \Upsilon_{\infty , t} (z_0)     \\ [10pt]
\Upsilon_{0, t}^\prime (z_0)                  &  \Upsilon_{\infty , t}^\prime (z_0)
\end{array} \right ]^{-1} 
\left [ \begin{array}{ c }
\widetilde \Psi_{k, t } (z_0)                       \\ [10pt]
\widetilde \Psi_{k, t }^\prime (z_0)                          
\end{array} \right ]
$$ 
This relation, and the $C^1$ uniform continuity in compact sets, of the state space and parameter space, of the
solutions of differential equations implies that the coefficients $d_0, d_\infty$ are continuous. 

Suppose that $\widetilde \Psi_{k, \hat t\ }(z_2) \not = 0$. Then $d_\infty(\ \hat t\ ) \not = 0$ and certainly  
$\Upsilon_{0, \hat t\ } (z_2) = 0$. This and the continuity of the situation imply that there exist
$\nu_1, \nu_2 > 0$ and a neighbourhood $U \subseteq [0,1]$ of $\hat t$ such that
$\nu_1 < \mid \widetilde \Psi_{k, t }(z) \mid$ for all $t \in U$ and all $z \in [z_2 - 2 \nu_2, z_2>$.

In particular, all the zeroes of $\widetilde \Psi_{k, \hat t }$ are in $[z_1, z_2 - \nu_2]$. 
Since it is a nontrivial solution of a second  order differential equation, it is true that 
$(\forall z)( z \in <z_1, z_2 - \nu_2]\ \wedge\  \widetilde \Psi_{k, \hat t }(z) = 0 
                     \Rightarrow \widetilde \Psi^\prime_{k, \hat t } (z) \not = 0)$.
This implies that, for $\mid t - \hat t \mid$ sufficiently small, all functions $\widetilde \Psi_{k, t }$ have 
the same number of zeroes (see item (9) above); contradicting the assumption that $\hat t$ was a point of
discontinuity for $ t \mapsto \sharp [(\widetilde \Psi_{k, t })^{-1}(\{ 0 \})]$. Hence 
$\widetilde \Psi_{k, \hat t }(z_2) = 0$. 

By applying the above constructions to all segments in a given order, the Proposition is proved     
\qed

\vspace*{2 true mm} 

\noindent {\bf Remarks.}

\vspace*{2 true mm} 
\noindent(1) It is not difficult to verify that, for $k \geq 2$, infinitesimal (nonzero) deformations obtained 
just by solving system $(\dagger \dagger \dagger \dagger )$ are nontrivial. Also, the ``Fourier'' nature of the
 solutions proposed says that nonzero solutions for different $k^\prime s$ are linearly independent.

\vspace*{2 true mm} 
\noindent(2) The ``glueing'' of function $f_1$ with $f_0$ in item (5) of the above Proposition may be much improved 
so as to make the subsequent arguments as smooth as pleased. But suppose $s_1 < s_2 < \cdots < s_{\ell +1}$, 
where $s_1 = s_{\ell +1}\ ({\rm mod}\ 2 \pi\ {\mathbb Z} )$, is a sequence of numbers corresponding to the points in 
${\mathcal C}(\zeta)$. The continuation process through the various $s_j$ $^\prime s$ can be done as smoothly as 
the solutions to system $(\clubsuit )$ are; except at the end when one finally reaches point $s_1$ again. 
No matter how smooth are the original data to the problem, in general, the glued solutions at $s_1$ are to be 
only of class $C^2$ (if the family $\{ \widetilde R_t \}_{t \in [0,1]}$ is $C^4$). 

\vspace*{2 true mm} 
\noindent(3) Any topological embedding $S^1 \hookrightarrow {\mathbb R}^2$ can certainly be perturbed by arbitrarily
small isotopies resulting in a $C^\infty$ embedding. Nevertheless, ``tracking differentiabilities'' may always prove 
to be of use in later situations (see Rado \cite{Rad01}).     

\vspace*{2 true mm} 
\noindent(4) Let $({\mathcal K}, \partial {\mathcal K}) \subseteq 
(\{ (x,z) \in {\mathbb R}^2 \mid x \geq 0 \} , \{ 0 \} \times {\mathbb R}) $ be diffeomorphic 
to the unit interval $([0,1], \partial [0,1])$ in such a way that ${\mathcal K}$ intersects the vertical 
axis orthogonally.
Then the surface of revolution obtained by revolving  ${\mathcal K}$ is topologically a sphere subject to 
the methods of this work. Indeed, first consider the one-manifold ${\mathcal K} \cup {\rm M}({\mathcal K})$, where 
$ {\rm M} : (x, z) \mapsto (-x, z)$. The analog of Lemma 2 can be obtained without difficulty as well as the 
succeeding Proposition. At the set of poles $\partial {\mathcal K}$ is where the fact that the height function
is Morse even analytic, and not just Morse analytic, becomes handy. These are points of null abscissa hence the
indicial equation has roots $1 \pm k$ and therefore, the best is to start the continuation process from there. 
Also, is easy to imagine examples of noncompact surfaces of revolution for which our main result is true.    
\label{sec:3}


\end{document}